\documentclass[a4paper,oneside,10pt]{article}%
\usepackage{amsmath}
\usepackage{amsfonts}
\usepackage{amssymb}
\usepackage{graphicx}
\usepackage{color}
\usepackage[square,numbers,sort&compress]{natbib}%
\providecommand{\U}[1]{\protect\rule{.1in}{.1in}}

\newtheorem{theorem}{Theorem}[section]

\newtheorem{corollary}[theorem]{Corollary}

\newtheorem{lemma}[theorem]{Lemma}

\newtheorem{remark}[theorem]{Remark}

\numberwithin{equation}{section}

\begin{document}

\title{A robust Kalman-Bucy filtering problem}
\author{Shaolin Ji\thanks{Zhongtai Institute of Finance, Shandong University, Jinan,
Shandong 250100, PR China. jsl@sdu.edu.cn. This research is supported by
National Natural Science Foundation of China (No. 11571203), the Programme of
Introducing Talents of Discipline to Universities of China (No. B12023).}
\quad Chuiliu Kong\thanks{Corresponding author. Zhongtai Institute of Finance,
Shandong University, Jinan, Shandong 250100, PR China.
kclsdmath@mail.sdu.edu.cn.}\quad Chuanfeng Sun\thanks{School of Mathematical
Sciences, University of Jinan, Jinan, Shandong 250022, P.R. China.
sms\_suncf@ujn.edu.cn. This research is partially supported by the National
Natural Science Foundation of China (No. 11701214), the Natural Science
Foundation of Shandong Province (No. ZR2017BA032).}}
\date{}
\maketitle

\textbf{Abstract}. A generalized Kalman-Bucy model under model uncertainty and
a corresponding robust problem are studied in this paper. We find that this
robust problem is equivalent to an estimate problem under a sublinear
operator. By Girsanov transformation and the minimax theorem, we prove that
this problem can be reformulated as a classical Kalman-Bucy filtering problem
under a new probability measure. The equation which governs the optimal
estimator is obtained. Moreover, the optimal estimator can be decomposed into
the classical optimal estimator and a term related to model uncertainty.

{\textbf{Key words}. }Kalman-Bucy filter, model uncertainty, robust, minimum
mean square estimator, minimax theorem, sublinear operator,

\section{Introduction}

It is well-known that Kalman and Bucy \cite{Kalman} built the fundamental
results of the filtering problem for linear systems, which are the foundation
of modern filtering theory (see Bensoussan \cite{Bensoussan}, Liptser and
Shiryaev \cite{Liptser} et al). Therefore, the filtering theory can be applied
to study stochastic optimal control problems with partial information (or
observation) in different fields. For example, Tang \cite{Tang} and
{\O }ksendal and Sulem \cite{Oksendal-Sulem} gave respectively the maximum
principle for partially observed optimal control problems\ of stochastic
differential equations and of forward-backward stochastic differential games;
Duncan and Pasik-Dunan \cite{Duncan1} and \cite{Duncan2} considered
respectively the optimal control for a partially observed linear stochastic
system with an exponential quadratic cost and with fractional brownian
motions; Bensoussan and Keppo \cite{Bensoussan-Keppo} and Lakner \cite{L}
considered the utility maximization problem under partial information where
the investor is interested in maximizing his utilities from consumption and
terminal wealth, and so on. Some fundamental research about forward-backward
stochastic differential equations can be referred to Ma, Protter and Yong
\cite{Ma-Protter-Yong} and Zhang \cite{Zhang}.

In more details, Kalman and Bucy considered that the signal process $({x}%
_{t})\in\mathbb{R}^{n}$ and the observation process $({m}_{t})\in
\mathbb{R}^{m}$ satisfy the following linear system:
\begin{equation}
\left\{
\begin{array}
[c]{rl}%
dx_{t} & =(F_{t}x_{t}+f_{t})dt+dw_{t},\\
x(0) & =x_{0},\\
dm_{t} & =(G_{t}{x}_{t}+g_{t})dt+dv_{t},\\
m(0) & =0
\end{array}
\right.  \label{eq-intro-system-1}%
\end{equation}
on a complete probability space $(\Omega,\mathcal{F},P)$. For the given
observation information $\mathcal{Z}_{t}=\sigma\{m(s),\ 0\leq s\leq t\}$, the
optimal estimator $\bar{x}_{t}$ of the signal ${x}_{t}$ solves the minimum
mean square estimation
\[
\min_{\zeta\in L_{\mathcal{Z}_{t}}^{2}(\Omega,P)}E_{P}\Vert x_{t}-\zeta
\Vert^{2}%
\]
where $L_{\mathcal{Z}_{t}}^{2}(\Omega,P)$ is the set of all the square
integral $\mathcal{Z}_{t}$-measurable random variables.

In this paper, we suppose that there exists model uncertainty for the system
(\ref{eq-intro-system-1}). Specifically, we don't know the true probability
$P$ and only know that it falls in a set of probability measures $\mathcal{P}%
$. For this case, it is naturally to consider the worst-case minimum mean
square estimation:%
\begin{equation}
\min_{\zeta}\sup_{P\in\mathcal{P}}E_{P}\Vert{x}_{t}-\zeta\Vert^{2}
\label{intro-problem}%
\end{equation}
which is to minimize the maximum expected loss over a range of possible
models, an idea that goes back at least as far as Wald \cite{Wald}. Recently,
this type of estimator has been utilized by Borisov \cite{Borisov1} and
\cite{Borisov2}, who studied the filtering of finite state Markov processes
with uncertainty of the transition intensity and the observation matrices.
Allan, Cohen \cite{Allan-Cohen} studied a parameter uncertainty problem in the
Kalman-Bucy filter with a control approach.

In our context, we adopt the $k$-ignorance model in Chen and Epstein
\cite{Chen-Epstein} to formulate the model uncertainty (readers can refer to
Section 2 for more details). Under this formulation, $\underset{P\in
\mathcal{P}}{\sup}E[\cdot]$ is a sublinear operator which is denoted by
$\mathcal{E}(\cdot)$. Thus, the problem (\ref{intro-problem}) can be
reformulated as
\[
\min_{\zeta}\mathcal{E(}\Vert{x}_{t}-\zeta\Vert^{2}).
\]
The main results about the\ estimation problem under sublinear operators are
obtained in Sun, Ji \cite{Sun-Ji} and Ji, Kong, Sun \cite{Ji-Kong-Sun}. Sun
and Ji \cite{Sun-Ji} studied this estimation problem for bounded random
variables. Ji, Kong and Sun \cite{Ji-Kong-Sun} deleted the boundedness
assumption and generalized the corresponding results to the case in which the
random variables belong to the space $L_{\mathcal{F}}^{2+\epsilon}(\Omega,P)$
where $\epsilon$ is a constant such that $\epsilon\in(0,1)$. The $k$-ignorance
model is one of the so-called drift-uncertainty models (see \cite{EJ-1,EJ-2}
for more general uncertainty models).

Under some mild conditions, we prove that the optimal estimator $\hat{x}$\ and
the optimal probability measure $P^{\theta^{\ast}}$ exist. It results that we
only need to consider the classical Kalman-Bucy filtering problem under the
probability measure $P^{\theta^{\ast}}$. Moreover, in some special cases, the
optimal estimator $\hat{x}_{t}$ can be decomposed to two parts. One part is
the optimal estimator of the signal process under the probability measure $P$
and the other part contains a parameter $\theta^{\ast}$ (see Corollary
\ref{decomposition} for details).

The paper is organized as follows. In section 2, a generalized robust
Kalman-Bucy filtering problem is introduced. The main results are given in
section 3. In section 4, we list some auxiliary theorems.

\section{Problem formulation}

Let $(\Omega,\mathcal{F},P)$ be a complete probability space on which two
independent $n$-dimensional and $m$-dimensional independent Brownian motions
$w(\cdot)$ and $v(\cdot)$ are defined. Assume that $\mathbb{F=}\{\mathcal{F}%
_{t},0\leq t\leq T\}$ is the $P$-augmentation of the natural filtration of
$w(\cdot)$ and $v(\cdot)$, where $\mathcal{F}={\mathcal{F}_{T}}$ and
$\mathcal{F}_{0}$ contains all $P$-null sets of $\mathcal{F}$. The means of
$w(\cdot)$ and $v(\cdot)$ are zero and the covariance matrices are $Q(\cdot)$
and $R(\cdot)$ respectively. The matrix $R(\cdot)$ is uniformly positive
definite. Denote by $\mathbb{R}^{n}$ the $n$-dimensional real Euclidean space
and $\mathbb{R}^{n\times k}$ the set of $n\times k$ real matrices. Let
$\langle\cdot,\cdot\rangle$ (resp. $\Vert\cdot\Vert$) denote the usual scalar
product (resp. usual norm) of $\mathbb{R}^{n}$ and $\mathbb{R}^{n\times k}$.
The scalar product (resp. norm) of $M=(m_{ij})$, $N=(n_{ij})\in\mathbb{R}%
^{n\times k}$ is denoted by $\langle M,N\rangle=tr\{MN^{\intercal}\}$ (resp.
$\Vert M\Vert=\sqrt{\langle M,M\rangle}$), where the superscript $^{\intercal
}$ denotes the transpose of vectors or matrices. For a $\mathbb{R}^{n}$-valued
vector $x=(x_{1},\cdot\cdot\cdot,x_{n})^{\intercal}$, $|x|:=(|x_{1}%
|,\cdot\cdot\cdot,|x_{n}|)^{\intercal}$; for two $\mathbb{R}^{n}$-valued
vectors $x$ and $y$, $x\leq y$ means that $x_{i}\leq y_{i}$ for $i=1,\cdot
\cdot\cdot,n$. Denote by $L_{\mathbb{F}}^{2}(0,T;\mathbb{R}^{n})$ the space of
$\mathbb{F}$-adapted $\mathbb{R}^{n}$-valued stochastic processes on $[0,T]$
such that
\[
\mathbb{E}_{P}\left[  \int_{0}^{T}|f(r)|^{2}dr\right]  <\infty.
\]
Through out this paper, $0$ denotes the matrix/vector with appropriate
dimension whose all entries are zero.

Under the probability measure $P$, the signal process $({x}_{t})\in
L_{\mathbb{F}}^{2}(0,T;\mathbb{R}^{n})$ and the observation process $({m}%
_{t})\in L_{\mathbb{F}}^{2}(0,T;\mathbb{R}^{m})$ satisfy
\begin{equation}
\left\{
\begin{array}
[c]{rl}%
dx_{t} & =(F_{t}x_{t}+f_{t})dt+dw_{t},\\
x(0) & =x_{0},\\
dm_{t} & =(G_{t}{x}_{t}+g_{t})dt+dv_{t},\\
m(0) & =0
\end{array}
\right.  \label{MO}%
\end{equation}
where $F_{t}\in\mathbb{R}^{n\times n},\ G_{t}\in\mathbb{R}^{m\times n}%
,\ f_{t}\in\mathbb{R}^{n},\ g_{t}\in\mathbb{R}^{m}$ are bounded, continuous
function in $t$ and $x_{0}\in\mathbb{R}^{n}$ be a given constant vector. Set
$\mathcal{Z}_{t}=\sigma\{m(s);0\leq s\leq t\}$. Then the filtration
$\mathbb{Z}\mathcal{=}\{\mathcal{Z}_{t},0\leq t\leq T\}$ represents the
observable information. By the Kalman-Bucy filtering theory (see Bensoussan
\cite{Bensoussan}, Kalman, Bucy \cite{Kalman}, Liptser and Shiryaev
\cite{Liptser}), the optimal estimate of $x_{t}$ under probability measure $P$
is governed by
\begin{equation}
\left\{
\begin{array}
[c]{ll}%
d\bar{x}_{t} & =(F_{t}\bar{x}_{t}+f_{t})dt+P_{t}G_{t}^{\top}R_{t}^{-1}%
dI_{t},\\
\bar{x}(0) & =x_{0},
\end{array}
\right.  \label{classical Kalman}%
\end{equation}
and the variance of estimate error $P_{t}=E_{P}[(x_{t}-\bar{x}_{t})(x_{t}%
-\bar{x}_{t})^{\top}]$ is governed by
\begin{equation}
\left\{
\begin{array}
[c]{rl}
& \frac{dP_{t}}{dt}=F_{t}P_{t}+P_{t}F_{t}^{\top}-P_{t}G_{t}^{\top}R_{t}%
^{-1}G_{t}P_{t}+Q_{t},\\
& P(0)=0
\end{array}
\right.  \label{Riccati}%
\end{equation}
where $\bar{x}_{t}=E_{P}(x_{t}|\mathcal{Z}_{t})$ and $I_{t}=m_{t}-\int_{0}%
^{t}(G_{s}\bar{x}_{s}+g_{s})ds$ is the so called innovation process which is a
Brownian motion adapted to $\mathbb{Z}$. Set $\mathcal{I}_{t}=\sigma
\{I(s);0\leq s\leq t\}$. Then the filtration $\{\mathcal{I}_{t}\}_{\{0\leq
t\leq T\}}$ equals to $\mathbb{Z}$.

Now we give the $k$-ignornace model which is proposed by Chen and Epstein
\cite{Chen-Epstein}. For a fixed $\mathbb{R}^{n}$-valued nonnegative constant
vector $\mu$, denote by $\Theta$ the set of all $\mathbb{R}^{n}$-valued
progressively measurable processes $(\theta_{t})$ with $|\theta_{t}|\leq\mu$.
Define%
\begin{equation}
\mathcal{P}=\{P^{\theta}\big|\frac{dP^{\theta}}{dP}=f_{T}^{\theta
}\ \;\text{for}\ \theta\in\Theta\} \label{CD}%
\end{equation}
where
\[
f_{T}^{\theta}=\exp\big(\int_{0}^{T}\theta_{t}^{\top}dw_{t}-\frac{1}{2}%
\int_{0}^{T}\Vert\theta_{t}\Vert^{2}dt\big).
\]
Due to the boundness of $\theta$, the Novikov's condition holds (see Karatzas
and Shreve \cite{K-S}). Therefore, $P^{\theta}$ defined by \eqref{CD} is a
probability measure and the processes $(w_{t}^{\theta})$ and $(v_{t})$ are
Brownian motions under this probability measure $P^{\theta}$ by Girsanov
theorem. Moreover, the probability measure $P^{\theta}$ is equivalent to the
probability measure $P$ with the Radon Nikodym derivative $\exp(\int_{0}%
^{T}\theta_{s}dw_{s}-\frac{1}{2}\int_{0}^{T}\theta_{s}^{2}ds)$. The
$k$-ignornace model describes an agent who is uncertain about the drift of the
underlying Brownian motion and allows any drift between $-\mu$ and $\mu$.

Taking into account the $k$-ignornace model, we generalize the Kalman-Bucy
filtering problem \eqref{MO} to the following minimax problem. Under every
probability measure $P^{\theta}\in\mathcal{P}$, consider
\begin{equation}
\left\{
\begin{array}
[c]{rl}%
dx_{t} & =(F_{t}x_{t}+f_{t}+\theta_{t})dt+dw_{t}^{\theta},\\
x(0) & =x_{0},\\
dm_{t} & =(G_{t}{x}_{t}+g_{t})dt+dv_{t},\\
m(0) & =0
\end{array}
\right.  \label{generalize KB}%
\end{equation}
where $w_{t}^{\theta}=w_{t}-\int_{0}^{t}\theta_{s}ds$ and study the minmax
problem
\begin{equation}
\inf_{\zeta\in L_{\mathcal{Z}_{t}}^{2+\epsilon}(\Omega,P,\mathbb{R}^{n})}%
\sup_{P^{\theta}\in\mathcal{P}}E_{P^{\theta}}\Vert{x}_{t}-\zeta\Vert^{2}
\label{robust problem}%
\end{equation}
where $\epsilon$ is a constant such that $0<\epsilon<1$ and $L_{\mathcal{Z}%
_{t}}^{2+\epsilon}(\Omega,P,\mathbb{R}^{n})$ is the set of all the
$\mathbb{R}^{n}$-valued $(2+\epsilon)$ integrable $\mathcal{Z}_{t}$-measurable
random variables.

It is easy to verify that $\mathcal{E}(\cdot)=\underset{P^{\theta}%
\in\mathcal{P}}{\sup}E_{P^{\theta}}[\cdot]$ is a sublinear operator. Thus, the
problem (\ref{robust problem}) can be represented as following: given the
observation information $\{\mathcal{Z}_{t}\}$, we want to find the optimal
estimator $\hat{x}_{t}$ for the signal $x_{t}$ for $t\in\lbrack0,T]$ such
that
\begin{equation}
\mathcal{E}\Vert{x}_{t}-\hat{x}_{t}\Vert^{2}=\inf_{\zeta\in\mathcal{K}_{t}%
}\mathcal{E}\Vert{x}_{t}-\zeta\Vert^{2} \label{OP}%
\end{equation}
where
\[
\mathcal{K}_{t}=\{\zeta:\Omega\rightarrow\mathbb{R}^{n};\ \zeta\in
L_{\mathcal{Z}_{t}}^{2+\epsilon}(\Omega,P,\mathbb{R}^{n})\}.
\]

Ji, Kong and Sun \cite{Ji-Kong-Sun} has explored the minimum mean square
estimator of random variables under sublinear operators and obtained the
existence and uniqueness results of the optimal estimator. In the next
section, we will utilize the results in \cite{Ji-Kong-Sun} to solve the
problem (\ref{OP}).

\begin{remark}
The optimal solution $\hat{x}_{t}$ of problem \eqref{OP} is called minimum
mean square estimator. It is also regarded as a minimax estimator in
statistical decision theory. If $\theta\equiv0$, then $\mathcal{P}$ contains
only the probability measure $P$ and the sublinear operator $\mathcal{E}%
(\cdot)$ degenerates to linear expectation operator. In this case, it is
well-known that the minimum mean square estimator $\hat{x}_{t}$ is just the
conditional expectation $E_{P}({x}_{t}|\mathcal{Z}_{t})$.
\end{remark}

\section{Main results}

In this section, we calculate the minimum mean square estimator $\hat{{x}}%
_{t}$ of the problem \eqref{OP} for $t\in\lbrack0,T]$. Without loss of
generality, all the statements in this section are only proved for the one
dimensional case.

\begin{lemma}
\label{LW} The set $\{\frac{dP^{\theta}}{dP}:P^{\theta}\in\mathcal{P}\}\subset
L^{1+\frac{2}{\epsilon}}(\Omega,\mathcal{F},P)$ is $\sigma(L^{1+\frac
{2}{\epsilon}}(\Omega,\mathcal{F},P),L^{1+\frac{\epsilon}{2}}(\Omega
,\mathcal{F},P))$-compact and the set $\mathcal{P}$ is convex.
\end{lemma}

\textbf{Proof}. Since $\theta$ is bounded, by Theorem \ref{AP} in Appendix,
the set $\{\frac{dP^{\theta}}{dP}:P^{\theta}\in\mathcal{P}\}$ is uniformly
normed bounded in $L^{1+\frac{2}{\epsilon}}(\Omega,\mathcal{F},P)$. From the
Theorem 4.1 of Chapter 1 in Simons \cite{Simons}, we know that the set
$\{\frac{dP^{\theta}}{dP}:P^{\theta}\in\mathcal{P}\}$ is $\sigma(L^{1+\frac
{2}{\epsilon}}(\Omega,\mathcal{F},P),L^{1+\frac{\epsilon}{2}}(\Omega
,\mathcal{F},P))$-compact.

The convexity of $\mathcal{P}$ is proved in Chen and Epstein
\cite{Chen-Epstein} (see Theorem 2.1 (c)).

\rightline{$\square$}

By Lemma \ref{LW}, we can apply the minimax theorem (see Theorem \ref{minmax})
to (\ref{robust problem}) which leads to the following theorem.

\begin{theorem}
\label{ER} For a given $t\in\lbrack0,T]$, there exists a $\theta^{\ast}%
\in\Theta$ such that
\[
\label{ERE}\inf_{\zeta\in\mathcal{K}_{t}}\mathcal{E}\Vert{x}_{t}-\zeta
\Vert^{2}=\inf_{\zeta\in\mathcal{K}_{t}}\sup_{P^{\theta}\in\mathcal{P}%
}E_{P^{\theta}}\Vert{x}_{t}-\zeta\Vert^{2}=\inf_{\zeta\in\mathcal{K}_{t}%
}E_{P^{\theta^{\ast}}}\Vert{x}_{t}-\zeta\Vert^{2}.
\]

\end{theorem}

\textbf{Proof}. Choose a sequence $\{\theta_{n}\}$, $n=1,2,\cdot\cdot\cdot$
such that%
\[
\lim_{n\rightarrow\infty}\inf_{\zeta\in\mathcal{K}_{t}}E_{P^{\theta_{n}}}%
[({x}_{t}-\zeta)^{2}]=\sup_{P^{\theta}\in\mathcal{P}}\inf_{\zeta\in
\mathcal{K}_{t}}E_{P^{\theta}}[({x}_{t}-\zeta)^{2}].
\]
Set $f_{n}=\frac{dP^{\theta_{n}}}{dP}$. By Koml\'{o}s theorem in the appendix
of Pham \cite{Pham}, we have that there exist a subsequence $\{f_{n_{k}%
}\}_{k\geq1}$ and a $f^{\ast}\in L^{1}(\Omega,\mathcal{F},P)$ such that
\[
\lim_{m\rightarrow\infty}\dfrac{1}{m}\sum_{k=1}^{m}f_{n_{k}}=f^{\ast
},\ P-a.s..
\]
Let $g_{m}=\dfrac{1}{m}\sum_{k=1}^{m}f_{n_{k}}$. We have $g_{m}%
\xrightarrow{P-a.s.}f^{\ast}$ and
\begin{equation}%
\begin{array}
[c]{rl}
& \sup_{P^{\theta}\in\mathcal{P}}\inf_{\zeta\in\mathcal{K}_{t}}E_{P^{\theta}%
}[({x}_{t}-\zeta)^{2}]=\lim_{n\rightarrow\infty}\inf_{\zeta\in\mathcal{K}_{t}%
}E_{P}[f_{n}({x}_{t}-\zeta)^{2}]=\lim_{k\rightarrow\infty}\inf_{\zeta
\in\mathcal{K}_{t}}E_{P}[f_{n_{k}}({x}_{t}-\zeta)^{2}]\\
& =\lim_{m\rightarrow\infty}\frac{1}{m}\sum_{k=1}^{m}\inf_{\zeta\in
\mathcal{K}_{t}}E_{P}[f_{n_{k}}({x}_{t}-\zeta)^{2}]\leq\lim_{m\rightarrow
\infty}\inf_{\zeta\in\mathcal{K}_{t}}\frac{1}{m}\sum_{k=1}^{m}E_{P}[f_{n_{k}%
}({x}_{t}-\zeta)^{2}]\\
& =\lim_{m\rightarrow\infty}\inf_{\zeta\in\mathcal{K}_{t}}E_{P}[g_{m}({x}%
_{t}-\zeta)^{2}].
\end{array}
\end{equation}

By Theorem \ref{AP}, for any given constants $p>1$ and $m$, we have
$\ E_{P}(g_{m})^{K}\leq M$ \ where $K=(1+\frac{2}{\epsilon})p$ and
$M=\exp(\frac{K^{2}-K}{2}\mu^{2}T)$. Then, we have $\left\{  |g_{m}%
|^{1+\frac{2}{\varepsilon}}:m=1,2,\cdot\cdot\cdot\right\}  $ is uniformly
integrable. Therefore, $\ g_{m}%
\xrightarrow{L^{1+\frac{2}{\epsilon}}(\Omega, \mathcal{F}, P)}f^{\ast}$ and
$f^{\ast}\in L^{1+\frac{2}{\epsilon}}(\Omega,\mathcal{F},P)$. According to the
convexity and weak compactness of the set $\{\frac{dP^{\theta}}{dP}:P^{\theta
}\in\mathcal{P}\}$, there exists a $\theta^{\ast}$ such that $\frac
{dP^{\theta^{\ast}}}{dP}=f^{\ast}$ and the following relations hold
\[%
\begin{array}
[c]{rl}%
\sup_{P^{\theta}\in\mathcal{P}}\inf_{\zeta\in\mathcal{K}_{t}}E_{P^{\theta}%
}[({x}_{t}-\zeta)^{2}] & \geq\inf_{\zeta\in\mathcal{K}_{t}}E_{P^{\theta^{\ast
}}}[({x}_{t}-\zeta)^{2}]\\
& =\inf_{\zeta\in\mathcal{K}_{t}}E_{P}[f^{\ast}({x}_{t}-\zeta)^{2}]\\
& =\inf_{\zeta\in\mathcal{K}_{t}}E_{P}[\lim_{m\rightarrow\infty}g_{m}({x}%
_{t}-\zeta)^{2}]\\
& \geq\limsup_{m\rightarrow\infty}\inf_{\zeta\in\mathcal{K}_{t}}E_{P}%
[g_{m}({x}_{t}-\zeta)^{2}]\\
& \geq\sup_{P^{\theta}\in\mathcal{P}}\inf_{\zeta\in\mathcal{K}_{t}%
}E_{P^{\theta}}[({x}_{t}-\zeta)^{2}]
\end{array}
\]
where the second $^{\prime}\geq^{\prime}$ is based on the upper
semi-continuous property. It follows that
\[
\sup_{P^{\theta}\in\mathcal{P}}\inf_{\zeta\in\mathcal{K}_{t}}E_{P^{\theta}%
}[({x}_{t}-\zeta)^{2}]=\inf_{\zeta\in\mathcal{K}_{t}}E_{P^{\theta^{\ast}}%
}[({x}_{t}-\zeta)^{2}].
\]
By the minimax theorem (Theorem \ref{minmax}), we obtain
\[
\sup_{P^{\theta}\in\mathcal{P}}\inf_{\zeta\in\mathcal{K}_{t}}E_{P^{\theta}%
}[({x}_{t}-\zeta)^{2}]=\inf_{\zeta\in\mathcal{K}_{t}}\sup_{P^{\theta}%
\in\mathcal{P}}E_{P^{\theta}}[({x}_{t}-\zeta)^{2}].
\]
It leads to that
\[
\inf_{\zeta\in\mathcal{K}_{t}}\sup_{P^{\theta}\in\mathcal{P}}E_{P^{\theta}%
}[({x}_{t}-\zeta)^{2}]=\inf_{\zeta\in\mathcal{K}_{t}}E_{P^{\theta^{\ast}}%
}[({x}_{t}-\zeta)^{2}].
\]
\rightline{$\square$}

Once we find the optimal $\theta^{\ast}$, the problem (\ref{generalize KB})
and (\ref{robust problem}) can be reformulated under the probability measure
$P^{\theta^{\ast}}$ by Theorem \ref{ER}. In more details, on the filtered
probability space $(\Omega,\mathcal{F},\{\mathcal{F}_{t}\}_{0\leq t\leq
T},P^{\theta^{\ast}})$, the processes $(x_{t})_{0\leq t\leq T}$ and
$(m_{t})_{0\leq t\leq T}$ satisfy
\begin{equation}
\left\{
\begin{array}
[c]{rl}%
{dx}_{t} & =(F_{t}{x}_{t}+f_{t}+\theta_{t}^{\ast})dt+dw_{t}^{\theta^{\ast}},\\
{x}(0) & =x_{0},\\
{dm}_{t} & =(G_{t}x_{t}+g_{t})dt+dv_{t},\\
{m}(0) & =0.
\end{array}
\right.  \label{MO1}%
\end{equation}
We solve the classical minimum mean square estimation problem under
$P^{\theta^{\ast}}$%
\begin{equation}
E_{P^{\theta^{\ast}}}\Vert{x}_{t}-\hat{x}_{t}\Vert^{2}=\inf_{\zeta
\in\mathcal{K}_{t}}E_{P^{\theta^{\ast}}}\Vert{x}_{t}-\zeta\Vert^{2}.
\label{NP2}%
\end{equation}

According to the above theorem, we study model \eqref{MO1} and the following
problem: to obtain the optimal estimator $\hat{x}_{t}$ such that
\begin{equation}
E_{P^{\theta^{\ast}}}\Vert{x}_{t}-\hat{x}_{t}\Vert^{2}=\inf_{\zeta\in
\bar{\mathcal{K}}_{t}}E_{P^{\theta^{\ast}}}\Vert{x}_{t}-\zeta\Vert^{2}
\label{KB}%
\end{equation}
where $\bar{\mathcal{K}}_{t}=\{\zeta:\Omega\rightarrow\mathbb{R}^{n}%
;\ \zeta\in L_{\mathcal{Z}_{t}}^{2}(\Omega,P^{\theta^{\ast}},\mathbb{R}%
^{n})\}$.

The problem \eqref{KB} is a classical linear partially observable system with
a fixed parameter $\theta^{\ast}$. This estimate problem is to characterize
the conditional distribution $P^{\theta^{\ast}}({x}_{t}\in A|\mathcal{Z}_{t}%
)$, where $A$ is a Borel set in $\mathbb{R}^{n}$. Then, we are in the realm of
Kalman-Bucy filtering and it is well known (Kalman and Bucy \cite{Kalman},
Liptser and Shiryaev \cite{Liptser}) that the conditional distribution is
again Gaussian and the conditional mean $\hat{x}_{t}=E_{P^{\theta^{\ast}}}%
({x}_{t}|\mathcal{Z}_{t})$ solves the following equation:%
\begin{equation}
\left\{
\begin{array}
[c]{rl}%
d\hat{x}_{t} & =(F_{t}\hat{x}_{t}+f_{t}+\widehat{\theta_{t}^{\ast}}%
)dt+P_{t}G_{t}R_{t}^{-1}d\hat{I}_{t},\\
\hat{x}_{t}(0) & =x_{0}%
\end{array}
\right.  \label{KBequation}%
\end{equation}
where $\widehat{\theta_{t}^{\ast}}=E_{P^{\theta^{\ast}}}[\theta^{\ast}%
_{t}|\mathcal{Z}_{t}]$, $\hat{I}_{t}={m}_{t}-\int_{0}^{t}(G_{s}\hat{x}%
_{s}+g_{s})ds$, $0\leq t\leq T$ is $\mathcal{Z}_{t}$-measurable Brownian
motion and the variance of the error equation $P_{t}=E_{P^{\theta^{\ast}}%
}[(x_{t}-\hat{x}_{t})^{2}|\mathcal{Z}_{t}]=E_{P^{\theta^{\ast}}}[(x_{t}%
-\hat{x}_{t})^{2}]$ satisfies the following equation.
\begin{equation}
\left\{
\begin{array}
[c]{rl}%
\frac{dP_{t}}{dt} & =F_{t}P_{t}+P_{t}F^{\intercal}_{t}+2E_{P^{\theta^{\ast}}%
}[\widehat{x_{t}\theta^{\ast\intercal}_{t}}-\hat{x}_{t}\widehat{\theta
^{\ast\intercal}_{t}}]-P_{t}G^{\intercal}_{t}R^{-1}_{t}G_{t}P_{t}+Q_{t},\\
P(0) & =0
\end{array}
\right.  \label{NewRiccati}%
\end{equation}
where $\widehat{x_{t}\theta^{\ast\intercal}_{t}}=E_{P^{\theta^{\ast}}}%
[x_{t}\theta^{\ast\intercal}_{t}|\mathcal{Z}_{t}]$. Thus, the optimal
estimator $\hat{x}_{t}$ of the problem \eqref{KB} has been obtained.

Next, we expound that this solution $\hat{x}_{t}$ is also the optimal
estimator of the problem \eqref{OP} at time $t\in\lbrack0,T]$.

\begin{theorem}
\label{ER1} Under the above assumptions, the solution $\hat{x}_{t}$ governed
by \eqref{KBequation} is also the optimal solution of the problem \eqref{OP}
at time $t\in\lbrack0,T]$.
\end{theorem}

\textbf{Proof.} Note that%
\begin{equation}
\inf_{\zeta\in\mathcal{K}_{t}}\sup_{P^{\theta}\in\mathcal{P}}E_{P^{\theta}%
}({x}_{t}-\zeta)^{2}=\inf_{\zeta\in\mathcal{K}_{t}}E_{P^{\theta^{\ast}}}%
({x}_{t}-\zeta)^{2}\geq\inf_{\zeta\in\bar{\mathcal{K}}_{t}}E_{P^{\theta^{\ast
}}}({x}_{t}-\zeta)^{2}\text{.}%
\end{equation}

Since $F_{t},\ G_{t},\ f_{t}\ \mbox{and}\ g_{t}$ are bounded continuous
functions in $t$ and $\theta^{\ast}$ is bounded, by Theorem \ref{ES}, $\hat
{x}_{t}$ is not only square integrable but also $(4+2\epsilon)$ integrable
under probability measure $P^{\theta^{\ast}}$. Then, the solution $\hat{x}%
_{t}$ to \eqref{KBequation} also belongs to $\mathcal{K}_{t}$. It yields that
$\hat{x}_{t}$\ is the optimal solution of problem \eqref{OP} at time
$t\in\lbrack0.T]$.\newline\rightline{$\square$}

\begin{corollary}
If the optimal $\theta^{\ast}_{t}$ is adapted to subfiltration $\mathcal{Z}%
_{t}$, then the optimal estimator $\hat{x}_{t}$ satisfies the following
simpler equation.
\begin{equation}
\left\{
\begin{array}
[c]{rl}%
d\hat{x}_{t} & =(F_{t}\hat{x}_{t}+f_{t}+\theta_{t}^{\ast})dt+P_{t}G_{t}%
R_{t}^{-1}d\hat{I}_{t},\\
\hat{x}_{t}(0) & =x_{0}%
\end{array}
\right.  \label{SKBequation}%
\end{equation}
where $P_{t}$ reduces to equation \eqref{Riccati}.
\end{corollary}

Define%
\[
R_{s}^{t}=P_{t}^{-1}[\Phi(t,s)Q(s)-\int_{s}^{t}A(t,r)G_{r}\Phi(r,s)Q(s)dr]
\]
where $A(t,s)=P_{s}G_{s}R_{s}^{-1}\exp^{\int_{s}^{t}(F_{r}-P_{r}G_{r}^{2}%
R_{r}^{-1})dr}$ is the impulse response of the classical Kalman-Bucy filter
and $\Phi(t,s)=\exp^{\int_{s}^{t}F_{r}dr}$. Then applying an analysis similar
to Theorem 1 in Davis and Varaiya \cite{Davis-Varaiya}, we obtain the
following Corollary.

\begin{corollary}
\label{decomposition} If the optimal $\theta^{\ast}_{t}$ is adapted to
subfiltration $\mathcal{Z}_{t}$, with equations \eqref{classical Kalman} and
\eqref{SKBequation}, then the optimal estimator $\hat{x}_{t}$ for any time
$t\in\lbrack0.T]$ can be expressed as
\begin{equation}
\hat{x}_{t}=\bar{x}_{t}+\int_{0}^{t}P_{t}R_{s}^{t}\theta_{s}^{\ast}ds.
\label{expression2}%
\end{equation}
where $\bar{x}_{t}$ is defined by equation \eqref{classical Kalman}.
\end{corollary}

\section{Appendix}

For the convenience of the reader, we put the following three theorems used in
the paper in this appendix.

\begin{theorem}
\textbf{(Girsanov \cite{Girsanov})} \label{AP} We suppose that $\phi
(t,\omega)$ satisfies the following conditions:\newline(1)\ $\phi(t,\omega)$
are measurable in both variables;\newline(2)\ $\phi(t,\omega)$ is
$\mathcal{F}_{t}$-measurable for fixed $t$;\newline(3)\ $\int_{0}^{T}%
|\phi(t,\omega)|^{2}dt<\infty$ almost everywhere; and $0<c_{1}\leq
|\phi(t,\omega)|\leq c_{2}$ for almost all $(t,\omega)$, then $\exp
[\alpha\zeta_{s}^{t}(\phi)]$ is integrable and for $\alpha>1$
\begin{equation}
\exp\big[\frac{(\alpha^{2}-\alpha)}{2}(t-s)c_{1}^{2}\big]\leq E[\exp
[\alpha\zeta_{s}^{t}(\phi)]]\leq\exp\big[\frac{(\alpha^{2}-\alpha)}%
{2}(t-s)c_{2}^{2}\big] \label{conclusion}%
\end{equation}
where $\zeta_{s}^{t}(\phi)=\int_{s}^{t}\phi(u,\omega)dw_{u}-\frac{1}{2}%
\int_{s}^{t}\phi^{2}(u,\omega)du$.
\end{theorem}

\begin{theorem}
[Pham \cite{Pham} Theorem B.1.2]\label{minmax} Let $\mathcal{X}$ be a convex
subset of a normed vector space, and $\mathcal{Y}$ be a convex subset of a
normed vector space $E$, compact for the weak topology $\sigma(E,E^{\prime})$.
Let $f:\mathcal{X}\times\mathcal{Y}\rightarrow R$ be a function
satisfying:\newline(1) $x\rightarrow f(x,y)$ is a continuous and convex on
$\mathcal{X}$ for all $y\in\mathcal{Y};$\newline(2) $y\rightarrow f(x,y)$ is a
concave on $\mathcal{Y}$ for all $x\in\mathcal{X}.$\newline Then, we have
\[
\sup_{y\in\mathcal{Y}}\inf_{x\in\mathcal{X}}f(x,y)=\inf_{x\in\mathcal{X}}%
\sup_{y\in\mathcal{Y}}f(x,y).
\]

\end{theorem}

\begin{theorem}
[Yong and Zhou \cite{Yong-Zhou}]\label{ES} Given a complete filtered
probability space $(\Omega,\mathcal{F},\{\mathcal{F}_{t}\}_{t\in\lbrack
0,T]},{P})$, on which a standard $\mathbb{R}$-valued Brownian motion
$W(\cdot)$ is defined. $\mathcal{F}_{t}=\sigma\{W_{s},s\leq t\}$. The process
$(X_{t})$ satisfies:
\[
\left\{
\begin{array}
[c]{rl}%
dX_{t} & =f(t,X_{t})dt+l(t,X_{t})dW_{t},\\
X_{0} & =x.
\end{array}
\right.
\]
Supposed that $f(t,x)$ and $l(t,x)$ satisfy the following conditions:\newline%
(L):\ $Lipschitz$ condition:\ $|f(t,x)-f(t,y)|+|l(t,x)-l(t,y)|\leq K|x-y|$,\ K
$\geq$ 0 is a constant;\newline(B):\ $\sup\limits_{t\in\lbrack0,T]}%
(|f(t,0)|+|l(t,0)|)<\infty$\newline and $E_{{P}}(|X_{0}|^{p})<\infty$,
$p\geq2$, then there exists constant $C_{p}$ such that
\[
E_{{P}}(\sup_{0\leq s\leq t}|X_{s}|^{p})\leq C_{p}(1+|X_{0}|^{p}%
)(1+t^{p})e^{C_{p}(t^{p/2}+t^{p})},\ t\geq0.
\]

\end{theorem}

\end{document}